\documentclass[11pt]{article}

\usepackage{amsfonts,amsmath,latexsym,color,epsfig,hyperref}
\newtheorem{theorem}{Theorem}
\newtheorem{lemma}{Lemma}

\newtheorem{problem}{Problem}

\def\qed{\ifvmode\mbox{ }\else\unskip\fi\hskip 1em plus 10fill$\Box$}

\input{epsf}
\makeatletter
\def\Ddots{\mathinner{\mkern1mu\raise\p@
\vbox{\kern7\p@\hbox{.}}\mkern2mu
\raise4\p@\hbox{.}\mkern2mu\raise7\p@\hbox{.}\mkern1mu}}
\makeatother

\title{\vspace{-0.7cm}Hedgehogs are not colour blind}
\author{David Conlon\thanks{Mathematical Institute, Oxford OX2 6GG, United Kingdom. Email: {\tt david.conlon@maths.ox.ac.uk}. Research supported by a Royal Society University Research Fellowship.} \and Jacob Fox\thanks{Department of Mathematics, Stanford University, Stanford, CA 94305, USA. Email: {\tt fox@math.mit.edu}. Research supported by a Packard Fellowship, by NSF Career Award DMS-1352121 and by an Alfred P. Sloan Fellowship.} \and Vojt\v ech R\"odl\thanks{Department of Mathematics and Computer Science, Emory University, Atlanta, GA 30322, USA. Email: {\tt rodl@mathcs.emory.edu}. Research partially supported by NSF grants DMS-1102086 and DMS-1301698.}}
\date{}

\begin{document}
\maketitle

\begin{abstract}
We exhibit a family of $3$-uniform hypergraphs with the property that their $2$-colour Ramsey numbers grow polynomially in the number of vertices, while their $4$-colour Ramsey numbers grow exponentially. This is the first example of a class of hypergraphs whose Ramsey numbers show a strong dependence on the number of colours.
\end{abstract}

\section{Introduction}

The Ramsey number $r_k(H)$ of a $k$-uniform hypergraph $H$ is the smallest $n$ such that any $2$-colouring of the edges of the complete $k$-uniform hypergraph $K_n^{(k)}$ contains a monochromatic copy of $H$. Similarly, for any $q \geq 2$, we may define a $q$-colour Ramsey number $r_k(H; q)$.

One of the main outstanding problems in Ramsey theory is to decide whether the Ramsey number for complete $3$-uniform hypergraphs is double exponential. The best known bounds, due to Erd\H{o}s, Hajnal and Rado~\cite{EHR65, ER52}, state that there are positive constants $c$ and $c'$ such that
\[2^{c t^2} \leq r_3(K_t^{(3)}) \leq 2^{2^{c' t}}.\]
Paul Erd\H{o}s has offered \$500 for a proof that the upper bound is correct, that is, that there exists a positive constant $c$ such that $r_3(K_t^{(3)}) \geq 2^{2^{c t}}$. Some evidence that this may be true was given by Erd\H{o}s and Hajnal (see, for example,~\cite{GRS90}), who showed that the analogous bound holds for $4$ colours, that is, that there exists a positive constant $c$ such that $r_3(K_t^{(3)}; 4) \geq 2^{2^{c t}}$. In this paper, we show that this evidence may not be so compelling by finding a natural class of hypergraphs, which we call hedgehogs, whose Ramsey numbers show a strong dependence on the number of colours. 

The hedgehog $H_t$ is the $3$-uniform hypergraph with vertex set $[t + \binom{t}{2}]$ such that for every pair $(i, j)$ with $1 \leq i < j \leq t$ there is a unique vertex $k > t$ such that $ijk$ is an edge. We will sometimes refer to the set $\{1, 2, \dots, t\}$ as the body of the hedgehog. Our main result is that the $2$-colour Ramsey number $r_3(H_t)$ grows as a polynomial in $t$, while the $4$-colour Ramsey number $r_3(H_t; 4)$ grows as an exponential in $t$. 

\begin{theorem} \label{thm:main} 
If $H_t$ is the $3$-uniform hedgehog with body of order $t$, then
\begin{itemize}
\item[(i)]
$r_3(H_t) \leq 4 t^3$,
\item[(ii)]
there exists a positive constant $c$ such that $r_3(H_t; 4) \geq 2^{ct}$.
\end{itemize}
\end{theorem}

For the intermediate $3$-colour case, we show that the answer is intimately connected with a special case of the multicolour Erd\H{o}s--Hajnal conjecture~\cite{EH89}. This conjecture states that for any complete graph $K$ with a fixed $q$-colouring of its edges, there exists a positive constant $c(K)$ such that any $q$-colouring (with the same $q$ colours) of the edges of the complete graph on $n$ vertices with no copy of $K$ contains a clique of order $n^{c(K)}$ which receives only $q-1$ colours. Though this conjecture is known to hold in a number of special cases (see, for example, Section~3.3 of~\cite{CFS15}), the best known general result, due to Erd\H{o}s and Hajnal themselves, says that there exists a positive constant $c(K)$ such that any $q$-colouring of the edges of the complete graph on $n$ vertices with no copy of $K$ contains a clique of order $e^{c(K) \sqrt{\log n}}$ which receives only $q-1$ colours.

We will be concerned with the particular case where $q = 4$ and the banned configuration $K$ is a rainbow triangle with one edge in each of the first three colours. 

\vspace{3mm}
{\bf Definition.} {\it Let $F(t)$ be the smallest $n$ such that every $4$-colouring of the edges of $K_n$, in red, blue, green and yellow, contains either a rainbow triangle $K$, with one edge in each of red, blue and green, or a clique of order $t$ with at most $3$ colours.} 

\vspace{3mm}
We will show that $r_3(H_t;3)$ is bounded above and below by polynomials in $F(t)$ (strictly speaking, the upper bound is a polynomial in $F(t^3)$, but, provided $F(t)$ does not jump pathologically, this will be at most polynomial in $F(t)$). Since the result of Erd\H{o}s and Hajnal mentioned in the previous paragraph implies that $F(t) \leq t^{c \log t}$ for some constant $c$, this in turn shows that $r_3(H_t; 3) \leq t^{c \log t}$ for some constant $c$. Moreover, the Erd\H{o}s--Hajnal conjecture holds in this case if and only if there is a polynomial upper bound for $r_3(H_t;3)$.

\begin{theorem} \label{thm:3col}
If $H_t$ is the $3$-uniform hedgehog with body of order $t$, then
\begin{itemize}
\item[(i)]
$r_3(H_t;3) =O(t^4 F(t^3)^2)$,
\item[(ii)]
$r_3(H_t; 3) \geq F(t)$.
\end{itemize}
In particular, there exists a constant $c$ such that $r_3(H_t;3) \leq t^{c \log t}$.
\end{theorem}

We will prove Theorem~\ref{thm:main} in the next section and Theorem~\ref{thm:3col} in Section~\ref{sec:3col}. We conclude by discussing a number of interesting questions that arose from our work.

\section{The basic dichotomy} \label{sec:basic}

In this section, we prove Theorem~\ref{thm:main}. We begin by proving that the $2$-colour Ramsey number of $H_t$ is at most $4 t^3$. 

\vspace{3mm}
{\bf Proof of Theorem~\ref{thm:main}(i):}
Let $n = 4t^3$. We will show that every red/blue-colouring of the complete $3$-uniform hypergraph on $n$ vertices contains a monochromatic copy of $H_t$. To begin, we define a partial colouring of the edges of the complete graph on the same vertex set. We will colour an edge $uv$ red if there are fewer than $\binom{t}{2} + t$ red triples containing $u$ and $v$. Similarly, we colour $uv$ blue if there are fewer than $\binom{t}{2} + t$ blue triples containing $u$ and $v$. To find a monochromatic $H_t$, it will clearly suffice to find a subset of order $t$ containing no red edge or no blue edge, since we can consider this set as the body of the hedgehog and embed the spines greedily.

 We claim that no vertex is contained in $2t^2$ red edges and $2t^2$ blue edges. Suppose, on the contrary, that $u$ is such a vertex and let $V_R$ and $V_B$ be the vertices which are connected to $u$ in red and blue, respectively. Since it is easy to see that no edge can be coloured both red and blue, $V_R$ and $V_B$ are disjoint. Moreover, for each vertex $v$ in $V_R$, since $uv$ is contained in fewer than $\binom{t}{2} + t$ red triples, there are at least 
\[|V_B| - \binom{t}{2} - t > \frac{|V_B|}{2}\] 
vertices $w$ in $V_B$ such that $uvw$ is blue. This implies that more than half of the triples $uvw$ with $v \in V_R$ and $w \in V_B$ are blue. However, by first considering vertices $w$ in $V_B$, the same argument also shows that more than half of these triples are red, a contradiction.

We now assign a colour to each vertex in the graph, colouring it red if it is contained in fewer than $2t^2$ red edges and blue otherwise. In the latter case, the claim of the last paragraph shows that it will be contained in fewer than $2t^2$ blue edges. By the pigeonhole principle, at least half the vertices in the graph have the same colour, say red. That is, we have a subset of order at least $n/2$ such that every vertex is contained in fewer than $2t^2$ red edges. By Brooks' theorem, we conclude that this set contains a subset of order $n/4t^2$ containing no red edge. Since $n/4t^2 \geq t$, this is the required set.
\qed

\vspace{3mm}
We will now show that the $4$-colour Ramsey number of $H_t$ is at least $2^{ct}$ for some positive constant $c$. This is clearly sharp up to the constant in the exponent.

\vspace{3mm}
{\bf Proof of Theorem~\ref{thm:main}(ii):}
A standard application of the first moment method gives a positive constant $c$ such that, for every integer $t \geq 4$, there is a $4$-colouring $\chi$ of the edges of the complete graph on $2^{ct}$ vertices with the property that every clique of order $t$ contains all $4$ colours. 

We now $4$-colour the edges of the complete $3$-uniform hypergraph on the same vertex set by colouring the triple $uvw$ with any colour which is not contained within the set $\{\chi(u, v), \chi(v,w), \chi(w,u)\}$. Suppose now that there is a monochromatic copy of $H_t$ with colour $1$, say, and let $u_1, u_2, \dots, u_t$ be the body of this copy. Then, in the original graph colouring $\chi$, none of the edges $u_i u_j$ with $1 \leq i < j \leq t$ received the colour $1$. However, this contradicts the property that every set of order $t$ contains all $4$ colours. 
\qed

\section{Three colours and the Erd\H{o}s--Hajnal conjecture} \label{sec:3col}

To prove Theorem~\ref{thm:3col}(i), we require two lemmas. The first is a result of Spencer~\cite{S72} which says that any $3$-uniform hypergraph with few edges contains a large independent set. 

\begin{lemma} \label{lem:Spencer}
If $H$ is a $3$-uniform hypergraph with $n$ vertices and $e$ edges, then $\alpha(H) = \Omega(n^{3/2}/e^{1/2})$.
\end{lemma}

The second lemma we require is a result of Fox, Grinshpun and Pach~\cite{FGP15} saying that the multicolour Erd\H{o}s--Hajnal conjecture holds for $3$-colourings of $K_n$ with no rainbow triangle. The result we use is somewhat weaker than the main result in~\cite{FGP15}, but will be more than sufficient for our purposes.

\begin{lemma} \label{lem:FGP}
Suppose that the edges of the complete graph $K_n$ have been $3$-coloured, in red, blue and green, so that there are no rainbow triangles with one edge in each of red, blue and green. Then there is a clique of order $n^{1/3}$ containing at most two of the three colours.
\end{lemma}

We are now ready to prove Theorem~\ref{thm:3col}(i), that $r_3(H_t;3) = O(t^4 F(t^3)^2)$.

\vspace{3mm}
{\bf Proof of Theorem~\ref{thm:3col}(i):}
Suppose that the edges of the complete $3$-uniform hypergraph on $n = c t^4 F(t^3)^2$ vertices have been $3$-coloured, in red, blue and green, where $c$ is a sufficiently large constant to be chosen later. We will $4$-colour the edges of the graph on the same vertex set as follows: if $u$ and $v$ are contained in fewer than $\binom{t}{2} + t$ triples of a given colour, then we give the edge $uv$ that colour, noting that an edge may receive more than one colour (but at most two). On the other hand, if an edge is not coloured with any of red, blue or green, we colour it yellow.

We claim that this colouring has at most $t^2 n^2$ triangles containing all three of the colours red, blue and green (where we include the possibility that two of these colours may appear on the same edge). To see this, note that there are at most $(\binom{t}{2} + t) \binom{n}{2}$ red triples containing a red edge. In particular, since the triangles we wish to count always contain a red edge, there are at most $(\binom{t}{2} + t) \binom{n}{2}$ of these triangles in the graph corresponding to a red triple. Since we may similarly bound the number of these triangles corresponding to blue or green triples, we see that, for $t \geq 3$, there are at most $3 (\binom{t}{2} + t) \binom{n}{2} \leq t^2 n^2$ triangles in the graph which contain all three of the colours red, blue and green, as required.

If we let $H$ be the $3$-uniform hypergraph on $n$ vertices whose edges correspond to triangles containg all three of the colours red, blue and green, Lemma~\ref{lem:Spencer} now yields a subset $U$ of order $\Omega(n^{1/2}/t )$ containing no such triangle. By taking $c$ to be sufficiently large, we may assume that $U$ has order at least $t F(t^3)$. 

We now consider the graph $G$ on vertex set $U$ whose edge set consists of all those edges which received two colours in the $4$-colouring defined above. If we fix a vertex $u \in U$, then each of the edges in $G$ that contain $u$ must have received the same two colours in the original colouring. Otherwise, we would have a triangle containing all three of the colours red, blue and green. Suppose, therefore, that every edge in $G$ that contains $u$ received the colours red and blue in the original colouring. Then, again using the property that every triangle contains at most two of the colours red, blue and green, we see that the neighbourhood of $u$ in $G$ contains no green edges. Therefore, if $u$ had $t$ neighbours in $G$, we could use this neighbourhood to find a green copy of $H_t$. Since a similar argument holds if the edges containing $u$ correspond to blue and green or to red and green, we may assume that every vertex $u \in U$ is contained in fewer than $t$ edges in the graph $G$. 

By Brooks' theorem, it follows that $U$ contains a subset $V$ of order at least $|U|/t \geq F(t^3)$ containing no edges from $G$, that is, such that every edge received at most one colour in the original colouring. Since $V$ is a $4$-coloured graph of order at least $F(t^3)$ containing no rainbow triangle in red, blue and green, there is a subset of order at least $t^3$ with at most three colours. If the missing colour is red, we may easily find a red copy of $H_t$ and similar conclusions hold if the missing colour is either blue or green. On the other hand, if the missing colour is yellow, we have a $3$-colouring, in red, blue and green, of a set of order at least $t^3$ containing no rainbow triangle, so Lemma~\ref{lem:FGP} tells us that there is a subset of order at least $t$ with at most two colours. If we again consider the missing colour, it is easy to find a monochromatic copy of $H_t$ in that colour.
\qed

\vspace{3mm}
The lower bound, $r_3(H_t; 3) \geq F(t)$ follows from a simple adaptation of the proof of Theorem~\ref{thm:main}(ii).

\vspace{3mm}
{\bf Proof of Theorem~\ref{thm:3col}(ii):}
By the definition of $F(t)$, there exists a $4$-colouring $\chi$, in red, blue, green and yellow, say, of the edges of the complete graph on $F(t) - 1$ vertices containing no rainbow triangle with one edge in each of red, blue and green and such that every clique of order $t$ contains all $4$ colours.

We now $3$-colour the complete $3$-uniform hypergraph on the same vertex set in red, blue and green, colouring the triple $uvw$ with any colour which is not contained within the set $\{\chi(u, v), \chi(v,w), \chi(w,u)\}$. Since there are no rainbow triangles in red, blue and green, this colouring is well-defined. Suppose now that there is a monochromatic copy of $H_t$ in red, say, and let $u_1, u_2, \dots, u_t$ be the body of this copy. Then, in the original graph colouring $\chi$, none of the edges $u_i u_j$ with $1 \leq i < j \leq t$ are red. However, this contradicts the property that every set of order $t$ contains all $4$ colours. 
\qed

\section{Concluding remarks}

The results of this paper raise a number of interesting questions, some of which we describe below.

\subsection{Higher uniformity hedgehogs}

The $k$-uniform hedgehog $H_t^{(k)}$ is the hypergraph with vertex set $[t + \binom{t}{k-1}]$ such that for every $(k-1)$-tuple $(i_1, \dots, i_{k-1})$ with $1 \leq i_1 <  \dots < i_{k-1} \leq t$ there is a unique vertex $i_k > t$ such that $i_1 \dots i_k$ is an edge. A straightforward generalisation of the proof of Theorem~\ref{thm:main}(i) gives the following result.

\begin{theorem} \label{thm:kuni}
For every integer $k \geq 4$, there exists a constant $c_k$ such that if $H_t^{(k)}$ is the $k$-uniform hedgehog with body of order $t$, then
\[r_k(H_t^{(k)}) \leq t_{k-2} (c_k t),\]
where the tower function $t_i(x)$ is defined by $t_1(x) = x$ and $t_{i+1}(x) = 2^{t_i(x)}$.
\end{theorem}

A construction due to Kostochka and R\"odl~\cite{KR06} shows that this result is tight for $k = 4$, that is, that there exists a positive constant $c$ such that $r_4(H_t^{(4)}) \geq 2^{ct}$. Since the construction is simple, we describe it in full. To begin, take a colouring of the edges of the complete graph on $2^{ct}$ vertices such that every set of order $t$ contains both a red triangle and a blue triangle. We then colour the edges of the $4$-uniform hypergraph on the same vertex set by colouring a $4$-tuple red if it contains a red triangle, blue if it contains a blue triangle and arbitrarily otherwise. It is easy to check that this $2$-colouring contains no monochromatic copy of $H_t^{(4)}$. Already for $k = 5$, we were unable to prove a matching lower bound, since it seems that one would first need to know how to prove a double-exponential lower bound for $r_3(K_t)$.

We were also unable to prove an analogue of Theorem~\ref{thm:main}(ii) for $k = 4$. Again, this is because of a basic gap in our understanding of hypergraph Ramsey problems. While we know that there are $4$-colourings of the $3$-uniform hypergraph on $2^{2^{ct}}$ vertices such that every subset of order $t$ receives at least two colours, we do not know if the following variant holds. 

\begin{problem}
Is there an integer $q$, a positive constant $c$ and a $q$-colouring of the $3$-uniform hypergraph on $2^{2^{ct}}$ vertices such that every subset of order $t$ receives at least three colours?
\end{problem}

A positive answer to the analogous question where we ask that every subset of order $t$ receives at least five colours would allow us to prove that there exists an integer $q$ such that $r_4(H_t^{(4)}; q) \geq 2^{2^{ct}}$. The proof of this statement is a variant of the proof of Theorem~\ref{thm:main}(ii). Indeed, suppose that we have a $q$-colouring $\chi$ of the edges of the $3$-uniform hypergraph $K_n^{(3)}$ such that every subset of order $t$ receives at least five colours. Then we define a colouring of the complete $4$-uniform hypergraph $K_n^{(4)}$ with at most $q + \binom{q}{2} + \binom{q}{3} + \binom{q}{4}$ colours by colouring the edge $uvwx$ with the set $\{\chi(uvw), \chi(vwx), \chi(wxu), \chi(xuv)\}$. It is now easy to check that if there is a monochromatic $H_t^{(4)}$ in this colouring, then, in the original colouring $\chi$, the body of the hedgehog is a subset of order $t$ which receives at most $4$ colours, contradicting our choice of $\chi$.

Our motivation for investigating higher uniformity hedgehogs was the hope that they might allow us to show that there are families of hypergraphs for which there is an even wider separation between the $2$-colour and $q$-colour Ramsey numbers. However, it seems likely that for hedgehogs the separation between the tower heights is at most one for any uniformity. This leaves the following problem open.

\begin{problem}
For any integer $h \geq 3$, do there exist integers $k$ and $q$ and a family of $k$-uniform hypergraphs for which the $2$-colour Ramsey number grows as a polynomial in the number of vertices, while the $q$-colour Ramsey number grows as a tower of height $h$?
\end{problem}

\subsection{Burr--Erd\H{o}s in hypergraphs}

The degeneracy of a graph $H$ is the minimum $d$ such that every induced subset contains a vertex of degree at most $d$. Building on work of Kostochka and Sudakov~\cite{KS03} and Fox and Sudakov~\cite{FS09}, Lee~\cite{L15} recently proved the famous Burr--Erd\H{o}s conjecture~\cite{BE75}, that graphs of bounded degeneracy have linear Ramsey numbers. That is, he showed that for every positive integer $d$ there exists a constant $c(d)$ such that the Ramsey number of any graph $H$ with $n$ vertices and degeneracy $d$ satisfies $r(H) \leq c(d) n$.

If we define the degeneracy of a hypergraph $H$ in a similar way, that is, as the minimum $d$ such that every induced subset contains a vertex of degree at most $d$, we may ask whether the analogous statement holds in hypergraphs. Unfortunately, as first observed by Kostochka and R\"odl~\cite{KR06}, the $4$-uniform analogue of the Burr--Erd\H{o}s conjecture is false, since $H_t^{(4)}$ is $1$-degenerate and $r_4(H_t^{(4)}) \geq 2^{ct}$. 

Since $H_t$ is a $1$-degenerate hypergraph, the results of this paper show that the Burr--Erd\H{o}s conjecture also fails for $3$-uniform hypergraphs and $3$ or more colours. For $4$ colours, this follows immediately from Theorem~\ref{thm:main}(ii). For $3$ colours, it follows from Theorem~\ref{thm:3col}(ii) and the observation that $F(t) = \Omega(t^3/\log^6 t)$. To show this, we amend a construction of Fox, Grinshpun and Pach~\cite{FGP15}, taking the lexicographic product of three $3$-colourings of the complete graph on $t/16 \log^2 t$ vertices, one for each triple of colours from the set $\{$red, blue, green, yellow$\}$ that contains yellow, each having the property that the union of any two colours contains no clique of order $4 \log t$. This colouring will contain no rainbow triangle with one edge in each of red, blue and green and no clique of order $t$ with at most $3$ colours. For further details, we refer the reader to Theorem~3.1 of~\cite{FGP15}.

While it is also unlikely that an analogue of the Burr--Erd\H{o}s conjecture holds in the $2$-colour case, it may still be the case that $r_3(H_t)$ is linear in the number of vertices, that is, that $r_3(H_t) = O(t^2)$. It would already be interesting to prove an approximate version of this statement.

\begin{problem}
Show that $r_3(H_t) = t^{2 + o(1)}$.
\end{problem}

\subsection{Multicolour Erd\H{o}s--Hajnal}

It is somewhat curious that our upper bound for $r_3(H_t; 3)$ mirrors the best known lower bound for $r_3(K_t;3)$, due to Conlon, Fox and Sudakov~\cite{CFS10}, which says that there exists a positive constant $c$ such that 
\[r_3(K_t;3) \geq 2^{t^{c \log t}}.\]
However, it seems likely that this is mere coincidence and that the function $F(t)$ defined in the introduction is polynomial in $t$. Phrasing the question in a more traditional fashion, we would very much like to know the answer to the following special case of the multicolour Erd\H{o}s--Hajnal conjecture.

\begin{problem}
Show that there exists a positive constant $c$ such that if the edges of $K_n$ are $4$-coloured, in red, blue, green and yellow, so that there are no rainbow triangles with one edge in each of red, blue and green, then there is a clique of order $n^c$ containing at most three of the four colours.

\end{problem}

That being said, if $F(t)$ were superpolynomial, it would not only disprove the multicolour Erd\H{o}s--Hajnal conjecture, it would also strengthen the curious correspondence between the bounds for $r_3(H_t; q)$ and $r_3(K_t; q)$. This would certainly be the more interesting outcome.

\end{document}